\documentclass[12pt,twoside]{amsart}
\usepackage{amsmath}
\usepackage{amsthm}
\usepackage{amsfonts}
\usepackage{amssymb}
\usepackage{latexsym}
\usepackage[all]{xy}
\usepackage{epsfig}
\usepackage{amscd,mathrsfs,graphicx,mathrsfs}

\date{}
\allowdisplaybreaks[4] \footskip=15pt

\renewcommand{\uppercasenonmath}[1]{}

\numberwithin{equation}{section} \theoremstyle{plain}
\newtheorem*{thm*}{Theorem A}
\newtheorem*{thm**}{Theorem B}
\newtheorem{thm}{Theorem}[section]
\newtheorem{cor}[thm]{Corollary}
\newtheorem*{cor*}{Corollary}
\newtheorem{lem}[thm]{Lemma}
\newtheorem*{lem*}{Lemma}
\newtheorem{prop}[thm]{Proposition}
\newtheorem*{prop*}{Proposition}

\newtheorem*{rem*}{Remark}

\newtheorem*{exa*}{Example}
\newtheorem{df}[thm]{Definition}
\newtheorem*{df*}{Definition}

\newtheorem*{conj*}{Conjecture}
\newtheorem*{ack*}{ACKNOWLEDGEMENTS}

\newcommand{\pf}{\noindent\begin {proof}}
\newcommand{\epf}{\end{proof}}

\begin{document}
\begin{center}
{\Large \bf Applications of cotorsion triples
 \footnotetext{
E-mail address: wren@cqnu.edu.cn.}}

\vspace{0.5cm}    Wei Ren\\
{\small School of Mathematical Sciences, Chongqing Normal University, Chongqing {\rm 401331}, China}
\end{center}


\bigskip
\centerline { \bf  Abstract}
\leftskip10truemm \rightskip10truemm
\noindent We study homotopy categories of model categories arising from a cotorsion triple, and the equivalences between corresponding stable categories. We characterize homological dimensions with respect to a cotorsion triple. Then, we lift cotorsion triple to complexes, and get the equivalence of homotopy categories of complexes via Quillen equivalence of model categories. Finally, we specify to Gorenstein cotorsion triple
$(\mathcal{GP}, \mathcal{W}, \mathcal{GI})$. By Quillen equivalence, it is shown that for a left-Gorenstein ring, there is an equivalence $\mathrm{K}(\mathcal{GP})\simeq \mathrm{K}(\mathcal{GI})$, which restricts to an equivalence $\mathrm{K}(\mathcal{P})\simeq \mathrm{K}(\mathcal{I})$ (compare to [J. Algebra, 2010, 324:2718-2731]); a new proof for Bennis and Mahdou's equality of Gorenstein global dimension is also given.
\bigskip


\leftskip0truemm \rightskip0truemm \vbox to 0.2cm{}

\section { \bf Introduction}
Throughout, let $\mathcal{A}$ be an abelian category with enough projectives and injectives. A pair of classes $(\mathcal{X}, \mathcal{Y})$ in $\mathcal{A}$ is a cotorsion pair provided that $\mathcal{X} =  {^\perp}\mathcal{Y}$ and $\mathcal{Y} = \mathcal{X}^{\perp}$, where $^{\perp}\mathcal{Y} = \{X \mid \mathrm{Ext}^{1}_{\mathcal{A}}(X, Y) = 0,~~\forall~~Y\in \mathcal{Y}\}$ and
$\mathcal{X}^{\perp} = \{Y \mid \mathrm{Ext}^{1}_{\mathcal{A}}(X, Y) = 0,~~\forall~~X\in \mathcal{X}\}$.
The cotorsion pair $(\mathcal{X}, \mathcal{Y})$ is said to be complete if for any $M\in \mathcal{A}$,
there exist short exact sequences $0\rightarrow Y\rightarrow X \rightarrow M \rightarrow 0$ and $0\rightarrow M\rightarrow Y' \rightarrow X' \rightarrow 0$ with $X, X'\in \mathcal{X}$ and $Y, Y'\in \mathcal{Y}$. Recall that a cotorsion pair $(\mathcal{X}, \mathcal{Y})$ is resolving if $\mathcal{X}$ is closed under taking kernels of epimorphisms between objects of $\mathcal{X}$, i.e. for any short exact sequence $0\rightarrow X^{'}\rightarrow X\rightarrow X^{''}\rightarrow 0$ with $X, X^{''}\in\mathcal{X}$, one has $X^{'}\in\mathcal{X}$. We say $(\mathcal{X}, \mathcal{Y})$ is coresolving if $\mathcal{Y}$ satisfies the dual statement. We say $(\mathcal{X}, \mathcal{Y})$ is hereditary if it is both resolving and coresolving. By \cite[Corollary 1.1.12]{Bec14}, a complete cotorsion pair is resolving if and only if it is coresolving.

A triple $(\mathcal{X}, \mathcal{Z}, \mathcal{Y})$ of classes of objects in $\mathcal{A}$ is called a (complete, hereditary) cotorsion triple, provided that both $(\mathcal{X}, \mathcal{Z})$ and $(\mathcal{Z}, \mathcal{Y})$ are (complete, hereditary) cotorsion pairs \cite{BR07, Chen10}.
For example, $(\mathcal{P}, \mathcal{A}, \mathcal{I})$ is a cotorsion triple, where $\mathcal{P}$ and $\mathcal{I}$ are respectively the subcategories of projective and injective objects in $\mathcal{A}$; if $R$ is a left-Gorenstein ring, then there is a cotorsion triple $(\mathcal{GP}, \mathcal{W}, \mathcal{GI})$ in $R$-Mod, where $\mathcal{GP}$ and $\mathcal{GI}$ are subcategories of Gorenstein projective and Gorenstein injective $R$-modules, and $\mathcal{W}$ consists of modules with finite projective dimension (equivalently, finite injective dimension). Beligiannis and Reiten \cite{BR07} studied torsion theory of stable categories by using cotorsion triples. Since $(\mathcal{GP}, \mathcal{GI})$ is a balanced pair (\cite[Definition 1.1]{Chen10}) arising from the Gorenstein cotorsion triple, the Gorenstein derived categories $\mathrm{D}_{\mathcal{GP}}(R)$ and $\mathrm{D}_{\mathcal{GI}}(R)$ are equivalent. Then in \cite{Chen10} Chen got a homotopy equivalence $\mathrm{K}(\mathcal{GP})\simeq \mathrm{K}(\mathcal{GI})$ by equivalence of natural composition functors $\mathrm{K}(\mathcal{GP})\stackrel{inc}\rightarrow \mathrm{K}(R)\rightarrow \mathrm{D}_{\mathcal{GP}}(R)$ and $\mathrm{K}(\mathcal{GI})\stackrel{inc}\rightarrow \mathrm{K}(R)\rightarrow \mathrm{D}_{\mathcal{GI}}(R)$. This restricts to an equivalence $\mathrm{K}(\mathcal{P})\simeq \mathrm{K}(\mathcal{I})$, which was firstly founded by Iyengar and Krause \cite{IK06} when the base ring is a commutative noetherian ring with dualizing complex.

In this paper, we study some new aspects of cotorsion triples. In Section 2, some properties of cotorsion triples are given. For example, in Proposition \ref{prop 3} we show that for a complete and hereditary cotorsion triple $(\mathcal{X}, \mathcal{Z}, \mathcal{Y})$, objects in $\mathcal{X}$ are ``projective'' with respect to $\mathcal{Y}$, i.e. for any object $M$, $M\in \mathcal{X}$ if and only if there is a completed commutative diagram
$$\xymatrix{&   &   & M\ar@{.>}[dl]_{\beta}\ar[d]^{\alpha} &   \\
 0 \ar[r]^{} &N \ar[r]^{f} &Y\ar[r]^{g} & L \ar[r]^{} &0,
  }$$
where $N$ is any object in $\mathcal{A}$ and $f: N\rightarrow Y$ is any left $\mathcal{Y}$-approximation of $N$. Dually, objects in $\mathcal{Y}$ are ``injective'' with respect to $\mathcal{X}$.

Recall that a model category, introduced by Quillen in \cite{Qui67}, is a category with three specified classes of morphisms,
called weak equivalences, cofibrations and fibrations, satisfying a few axioms \cite{DS95, Hov99}. For a bicomplete model category, the associated homotopy category is constructed by formally inverting the weak equivalences, i.e. localization with respect to weak equivalences. In Section 3, we study stable categories arising from cotorsion triples. By Hovey's correspondence between complete cotorsion pairs and model structures \cite[Theorem 2.2]{Hov02}, we show that there is a projective model structure and an injective model structure induced from a complete and hereditary cotorsion triple $(\mathcal{X}, \mathcal{Z}, \mathcal{Y})$, with the identity adjunction being a Quillen equivalence between them; and consequently, there are equivalences between associated homotopy categories and stable categories (see Theorem \ref{thm 1}):
$$ \mathcal{X}/\mathcal{P}\simeq \mathrm{Ho}(\mathcal{M}^{proj})\simeq \mathrm{Ho}(\mathcal{M}^{inj})\simeq \mathcal{Y}/\mathcal{I}.$$ Moreover, we characterize morphisms in the associated homotopy categories, and when the objects are equivalent in the stable categories; see Proposition
\ref{prop 6} and \ref{prop 7}.

In Section 4, relative homological dimensions with respect to the complete and hereditary cotorsion triple $(\mathcal{X}, \mathcal{Z}, \mathcal{Y})$ are investigated. It is immediate from Proposition \ref{prop 3} and \ref{prop 4} that $\mathrm{Hom}_{\mathcal{A}}(-, -)$ is right balanced
by $\mathcal{X}\times \mathcal{Y}$ in the sense of \cite[Definition 8.2.13]{EJ00}. Thus we can define relative derived functors $\mathrm{Ext}^{i}_{\mathcal{XY}}(-,-)$ to characterize $\mathcal{Z}$-projective and $\mathcal{Z}$-injective dimensions (Theorem \ref{thm 2} and \ref{thm 3}). Consequently, we get that global $\mathcal{Z}$-projective dimension and global $\mathcal{Z}$-injective dimension of $\mathcal{A}$ are equal (Corollary \ref{cor 1}).

Cotorsion triple of complexes and homotopy equivalences are studied in Section 5. If $\mathcal{A}$ is of finite global $\mathcal{Z}$-dimension, we can lift the complete and hereditary cotorsion triple $(\mathcal{X}, \mathcal{Z}, \mathcal{Y})$ in $\mathcal{A}$ to a complete and hereditary cotorsion triple $(dw\widetilde{\mathcal{X}}, \widehat{\mathcal{Z}}, dw\widetilde{\mathcal{Y}})$ in the category $\mathrm{Ch}(\mathcal{A})$ of complexes. Thus, there is a projective model structure and an injective model structure on $\mathrm{Ch}(\mathcal{A})$, and then we get the equivalence $\mathrm{K}(\mathcal{X})\simeq \mathrm{K}(\mathcal{Y})$ by Quillen equivalence of model categories; see Theorem \ref{thm 4}. We remark that the equivalence was proved by Chen \cite[Theorem A]{Chen10} via relative derived categories.

We specify to the Gorenstein cotorsion triple in Section 6. We note that Gorenstein cotorsion triple are obtained under different conditions, see for example \cite{Bel00, EJ00, Hov02, ZAD14}. If the ring is left-Gorenstein \cite{Bel00} or virtually Gorenstein \cite{ZAD14}, then it follows from Theorem \ref{thm 1} that there is an equivalence $\mathcal{GP}/\mathcal{P}\simeq \mathcal{GI}/\mathcal{I}$ between stable categories. By the results in Section 4, we can give a new proof for Bennis and Mahdou's Gorenstein global dimension \cite{BM10}, where they proved by so called strongly Gorenstein projective and injective modules (analogous to free module). Our argument, similar to the way of determining the global dimension in classical homological algebra, is based on vanishing of relative derived functors. This serves as an example to support the assertion of Holm \cite{Holm} that ``every result in classical homological algebra has a counterpart in Gorenstein homological algebra''. Moreover, we can lift the Gorenstein cotorsion triple to complexes, and then we get \cite[Theorem B]{Chen10} (Corollary \ref{cor 6}): Let $R$ be a left-Gorenstein ring, then there is an equivalence of categories $\mathrm{K}(\mathcal{GP})\simeq \mathrm{K}(\mathcal{GI})$, which restricts to an equivalence $\mathrm{K}(\mathcal{P})\simeq \mathrm{K}(\mathcal{I})$.
Moreover, consider the compact objects and the equivalence $\mathrm{K}(\mathcal{P})^{c}\simeq \mathrm{D}^{b}(R^{op}\text{-}\mathrm{mod})$ of Neeman \cite[Proposition 7.12]{Nee08} (compare to J{\o}rgensen \cite[Theorem 3.2]{Jor05}), and equivalence $\mathrm{K}(\mathcal{I})^{c}\simeq \mathrm{D}^{b}(R\text{-}\mathrm{mod})$ by Krause \cite[Proposition 2.3(2)]{Kra05}, then we have:
Let $R$ be a left-Gorenstein ring which is left noetherian and right coherent, then there is a duality $\mathrm{D}^{b}(R^{op}\text{-}\mathrm{mod})\simeq \mathrm{D}(R\text{-}\mathrm{mod})$ of triangulated categories; see Corollary \ref{cor 7} or \cite[Corollary C]{Chen10}.

\section { \bf Some properties of cotorsion triples}

In this section, some properties of cotorsion triple are presented. We start by recalling the following definition.

\begin{df} \label{df 1}(\cite{Chen10})
A triple $(\mathcal{X}, \mathcal{Z}, \mathcal{Y})$ of classes of objects in an abelian category $\mathcal{A}$ is called a cotorsion triple, provided that both $(\mathcal{X}, \mathcal{Z})$ and $(\mathcal{Z}, \mathcal{Y})$ are cotorsion pairs. The cotorsion triple is complete (resp. hereditary) provided that both of the two cotorsion pairs $(\mathcal{X}, \mathcal{Z})$ and $(\mathcal{Z}, \mathcal{Y})$ are complete (resp. hereditary).
\end{df}

Note that Beligiannis and Reiten \cite{BR07} also used the notion of cotorsion triple in studying torsion theories of triangulated categories, where it exactly means a complete cotorsion triple. We have the following observations on cotorsion triples.

Recall that a full subcategory $\mathcal{C}$ of $\mathcal{A}$ is thick if it is closed under direct summands and satisfies the 2-out-of-3 property:
for every exact sequence $0\rightarrow A\rightarrow B \rightarrow C \rightarrow 0$ in $\mathcal{A}$ with two terms in $\mathcal{C}$,
the third term belongs to $\mathcal{C}$ as well.

\begin{prop}\label{prop 1}
Let $(\mathcal{X}, \mathcal{Z}, \mathcal{Y})$ be a complete cotorsion triple. Then $\mathcal{Z}$ is thick if and only if $(\mathcal{X}, \mathcal{Z}, \mathcal{Y})$ is hereditary.
\end{prop}

\begin{proof}
It is easy to see that $\mathcal{Z}$ is thick if and only if the cotorsion pair $(\mathcal{X}, \mathcal{Z})$ is coresolving and $(\mathcal{Z}, \mathcal{Y})$ is resolving. By \cite[Corollary 1.1.12]{Bec14}, a complete cotorsion pair is resolving if and only if it is coresolving. This implies that $(\mathcal{X}, \mathcal{Z}, \mathcal{Y})$ is hereditary.
\end{proof}

\begin{prop}\label{prop 2}
Let $(\mathcal{X}, \mathcal{Z}, \mathcal{Y})$ be a complete and hereditary cotorsion triple. Then $\mathcal{X} \cap \mathcal{Z} = \mathcal{P}$, $\mathcal{Z} \cap \mathcal{Y} = \mathcal{I}$.
\end{prop}

\begin{proof}
Since  $(\mathcal{X}, \mathcal{Z})$ and $(\mathcal{Z}, \mathcal{Y})$ are cotorsion pairs, it is easy to see that $\mathcal{P} \subseteq \mathcal{X} \cap \mathcal{Z}$. Now assume that $M$ is any object in $\mathcal{X} \cap \mathcal{Z}$. Consider the exact sequence $0\rightarrow K\rightarrow P\rightarrow M\rightarrow 0$, where $P$ is projective. We have $K\in \mathcal{Z}$  by the thickness of $\mathcal{Z}$ since $M$ and $P$ are in $\mathcal{Z}$, and then $\mathrm{Ext}_{\mathcal{A}}^{1}(M, K) = 0$. This implies that the sequence is split, and hence $M$ is projective. Similarly, we can prove $\mathcal{Z} \cap \mathcal{Y} = \mathcal{I}$.
\end{proof}

In the following, we give characterizations for objects in the subcategories $\mathcal{X}$ and $\mathcal{Y}$ of a cotorsion triple
$(\mathcal{X}, \mathcal{Z}, \mathcal{Y})$ by commutative diagrams. Roughly speaking,
it is shown that any object in $\mathcal{X}$ is ``projective'' with respect to $\mathcal{Y}$
and any object in $\mathcal{Y}$ is ``injective'' with respect to $\mathcal{X}$.

Recall that a morphism $f: M\rightarrow Y$ (or briefly, $Y$) is called a left $\mathcal{Y}$-approximation (or $\mathcal{Y}$-preenvelope) of $M$,
if $Y\in \mathcal{Y}$ and any morphism from $M$ to an object in $\mathcal{Y}$ factors through $f$, i.e. for any
$Y^{'}\in \mathcal{Y}$ the morphism
$\mathrm{Hom}_{\mathcal{A}}(f, Y^{'}): \mathrm{Hom}_{\mathcal{A}}(Y, Y^{'})\rightarrow \mathrm{Hom}_{\mathcal{A}}(M, Y^{'})$ is surjective. In particular,  $f: M\rightarrow Y$ is said to be a special left $\mathcal{Y}$-approximation if $f$ is monic and
$\mathrm{Coker}(f)\in {^{\perp}}\mathcal{Y}$.
Dually, one has the notions of right $\mathcal{X}$-approximation (or $\mathcal{X}$-precover) and special right $\mathcal{X}$-approximation.

\begin{lem}\label{lem 1}
Let $(\mathcal{X}, \mathcal{Z}, \mathcal{Y})$ be a hereditary cotorsion triple and $M$, $N$ be objects in $\mathcal{A}$.
If $M\in \mathcal{X}$, then for any special left $\mathcal{Y}$-approximation $f: N\rightarrow Y$,
any homomorphism $\alpha: M\rightarrow L = \mathrm{Coker}(f)$
can be lifted to $\beta: M\rightarrow Y$. That is, we have the following completed commutative diagram:
$$\xymatrix{&   &   & M\ar@{.>}[dl]_{\beta}\ar[d]^{\alpha} &   \\
 0 \ar[r]^{} &N \ar[r]^{f} &Y\ar[r]^{g} & L \ar[r]^{} &0.
  }$$
\end{lem}

\begin{proof}
Note that in the exact sequence $0\rightarrow N \stackrel{f}\rightarrow Y \stackrel{g}\rightarrow L\rightarrow 0$, $L\in \mathcal{Z}$.
Consider the exact sequence $0\rightarrow K \stackrel{\iota}\rightarrow P \stackrel{\pi}\rightarrow L\rightarrow 0$, where $P$ is projective.
It is easy to see that there exist homomorphisms $\delta: P\rightarrow Y$ and $\gamma: K\rightarrow N$ such that the following diagram commutes:
$$\xymatrix{0 \ar[r]^{}  &K\ar[d]_{\gamma} \ar[r]^{\iota} & P \ar[d]_{\delta}
  \ar[r]^{\pi} &L\ar@{=}[d] \ar[r]^{ } &0 \\
  0 \ar[r]^{ } & N \ar[r]^{f} & Y \ar[r]^{g} & L \ar[r]^{ } & 0.
  }$$
Since $P$ and $L$ are in $\mathcal{Z}$, we have $K\in\mathcal{Z}$.
Hence $\mathrm{Ext}_{\mathcal{A}}^{1}(M, K) = 0$, and it yields that for any
homomorphism $\alpha: M\rightarrow L$, there exists a homomorphism $\sigma: M\rightarrow P$ such that $\alpha = \pi \sigma$. Put
$\beta = \delta\sigma$, then $\alpha = g\beta$, and this completes the proof.
\end{proof}

\begin{prop}\label{prop 3}
Let $(\mathcal{X}, \mathcal{Z}, \mathcal{Y})$ be a complete and hereditary cotorsion triple and $M$ an object in $\mathcal{A}$.
The following are equivalent:\\
\indent$\mathrm{(1)}$  $M \in \mathcal{X}$.\\
\indent$\mathrm{(2)}$  For any object $N\in \mathcal{A}$ and any left
$\mathcal{Y}$-approximation $f: N\rightarrow Y$, any homomorphism $\alpha: M\rightarrow L = \mathrm{Coker}(f)$
can be lifted to $\beta: M\rightarrow Y$. That is, we have the following completed commutative diagram:
$$\xymatrix{&   &   & M\ar@{.>}[dl]_{\beta}\ar[d]^{\alpha} &   \\
 0 \ar[r]^{} &N \ar[r]^{f} &Y\ar[r]^{g} & L \ar[r]^{} &0.
  }$$
\end{prop}

\begin{proof}
(1)$\Longrightarrow$(2). Since the cotorsion pair $(\mathcal{Z}, \mathcal{Y})$ is complete, any object $N\in \mathcal{A}$
has a special left $\mathcal{Y}$-approximation $f^{'}: N\rightarrow Y^{'}$.
Then we have the following commutative diagram:
$$\xymatrix{0 \ar[r]  &N\ar@{=}[d] \ar[r]^{f} & Y \ar[d]_{\gamma}\ar[r]^{g}
&L\ar[d]_{\delta} \ar[r]^{ } &0 \\
0 \ar[r]^{ } & N \ar[r]^{f'} & Y' \ar[r]^{g'} & L' \ar[r]^{} &
0}$$
For $\delta\alpha: M\rightarrow L^{'}$, there exists, by Lemma \ref{lem 1}, a homomorphism $\beta^{'}: M\rightarrow Y^{'}$
such that $\delta\alpha = g^{'}\beta^{'}$. Moreover, since the right square
is a pullback of $g^{'}$ and $\delta$, we have a homomorphism $\beta: M\rightarrow Y$ satisfying $\alpha = g \beta$ and
$\beta^{'} = \gamma \beta$.

(2)$\Longrightarrow$(1).  Let $Z$ be any object in $\mathcal{Z}$. It suffices to prove that any extension
$0\rightarrow Z \stackrel{f}\rightarrow N \stackrel{g}\rightarrow M\rightarrow 0$ of $Z$ by $M$ is split.
For $Z$, there is an exact sequence $0\rightarrow Z \stackrel{i}\rightarrow I \stackrel{p}\rightarrow L\rightarrow 0$
with $I$ injective. By Proposition \ref{prop 2},  $I\in \mathcal{Z}\cap \mathcal{Y}$. Moreover, we have
$L\in \mathcal{Z}$, and this implies that $i: Z\rightarrow I$ is a special left $\mathcal{Y}$-approximation.
By the injectivity of $I$, we have the following commutative diagram:
 $$\xymatrix{0 \ar[r]  &Z\ar@{=}[d] \ar[r]^{f} & N \ar[d]_{\gamma}\ar[r]^{g}
&M\ar[d]_{\delta} \ar[r]^{ } &0 \\
0 \ar[r]^{ } & Z \ar[r]^{i} & I \ar[r]^{p} & L \ar[r]^{} &
0}$$
For $\delta: M\rightarrow L$, there exists a homomorphism $\alpha: M\rightarrow I$ such that
$\delta = p \alpha$. Note that the right square is a pullback, and then there is a homomorphism $\beta: M\rightarrow N$
such that $\alpha = \gamma \beta$ and $g \beta = 1_M$. Hence the upper exact sequence is split as desired.
\end{proof}

Dually, we have the following.

\begin{prop}\label{prop 4}
Let $(\mathcal{X}, \mathcal{Z}, \mathcal{Y})$ be a complete and hereditary cotorsion triple and $M$ an object in  $\mathcal{A}$.
The following are equivalent:\\
\indent$\mathrm{(1)}$  $M\in \mathcal{Y}$.\\
\indent$\mathrm{(2)}$  For any object $N\in \mathcal{A}$ and any right
$\mathcal{X}$-approximation $f: X\rightarrow N$,
any homomorphism $\alpha: K = \mathrm{Ker}(f)\rightarrow M$
can be lifted to $\beta: X\rightarrow M$. That is, we have the following completed commutative diagram:
$$\xymatrix{&  M  \\
 0 \ar[r]^{} &K \ar[u]^{\alpha}\ar[r]^{g} &X\ar@{.>}[lu]_{\beta}\ar[r]^{f} & N \ar[r]^{} &0.
  }$$
\end{prop}

\section {\bf Stable categories arising from cotorsion triples}

A model structure on an abelian category $\mathcal{A}$ consists of three distinguished classes of maps, called weak equivalences, cofibrations and fibrations respectively, satisfying a few axioms. A model category is a bicomplete category (i.e. a category possessing arbitrary small limits and colimits) equipped with a model structure. For model categories we refer to the original source of Quillen \cite{Qui67}, as well as \cite{DS95, Hir03, Hov99}.

By Hovey's correspondence \cite[Theorem 2.2]{Hov02}, an abelian model structure on $\mathcal{A}$ is equivalent to a triple $[\mathcal{A}_{c}, \mathcal{A}_{tri}, \mathcal{A}_{f}]$ of subcategories, for which $\mathcal{A}_{tri}$ is thick and both $(\mathcal{A}_{c}, \mathcal{A}_{f}\cap \mathcal{A}_{tri})$ and $(\mathcal{A}_{c}\cap \mathcal{A}_{tri}, \mathcal{A}_{f})$ are complete cotorsion pairs. In this case, $\mathcal{A}_{c}$ is the class of cofibrant objects, $\mathcal{A}_{tri}$ is the class of trivial objects and $\mathcal{A}_{f}$ is the class of fibrant objects. The model structure is called ``abelian'' since it is compatible with the abelian structure of the category in the following way: (trivial) cofibrations are monomorphisms with (trivially) cofibrant cokernel, (trivial) fibrations are epimorphisms with (trivially) fibrant kernel, and weak equivalences are morphisms which factor as a trivial cofibratin followed by a trivial fibration.

For convenience, $[\mathcal{A}_{c}, \mathcal{A}_{tri}, \mathcal{A}_{f}]$ is called a Hovey triple. In particular, if $\mathcal{A}_{c} = \mathcal{A}$, i.e. every object in $\mathcal{A}$ is cofibrant, then it is called an injective Hovey triple. Dually, if $\mathcal{A}_{f} = \mathcal{A}$, then it is called a projective Hovey triple. We will not distinguish a Hovey triple and the corresponding model structure. Now we have the following model structures on $\mathcal{A}$ immediately by Hovey's correspondence.

\begin{prop}\label{prop 5}
Let $\mathcal{A}$ be an abelian category with a complete and hereditary cotorsion triple $(\mathcal{X}, \mathcal{Z}, \mathcal{Y})$. Then there are both a projective model structure $\mathcal{M}^{proj}=[\mathcal{X}, \mathcal{Z}, \mathcal{A}]$ and an injective model structure $\mathcal{M}^{inj} = [\mathcal{A}, \mathcal{Z}, \mathcal{Y}]$ on $\mathcal{A}$.
\end{prop}

For a bicomplete abelian category $\mathcal{A}$ with the model structure $\mathcal{M} = [\mathcal{A}_{c}, \mathcal{A}_{tri}, \mathcal{A}_{f}]$, the associated homotopy category $\mathrm{Ho}(\mathcal{A}, \mathcal{M})$ is constructed by formally inverting the weak equivalences, i.e. the localization with respect to weak equivalences. There is no risk of ambiguity that we simply use the notation $\mathrm{Ho}(\mathcal{M})$ to emphasize the model structure if the underlying category is apparent, and we use the notation $\mathrm{Ho}(\mathcal{A})$ when we need not to declare the model structure on $\mathcal{A}$.

By applying the factorization axiom (see \cite[Definition 1.1.3]{Hov99} or \cite[Definition 3.3, MC5]{DS95}) of the model
structure, one can factor the map $0\rightarrow M$ to obtain a trivial fibration $p_{M}: QM\rightarrow M$ with $QM$ cofibrant,
and  one can factor the map $M\rightarrow 0$ to obtain a trivial cofibration $i_{M}: M\rightarrow RM$ with $RM$ fibrant.
If $M$ is itself cofibrant, let $QM = M$; if $M$ is itself fibrant, let $RM = M$. $QM$ and $RM$ are called cofibrant and fibrant
replacement of $M$ respectively. Since $(\mathcal{A}_{c}, \mathcal{A}_{f}\cap \mathcal{A}_{tri})$ and $(\mathcal{A}_{c}\cap \mathcal{A}_{tri}, \mathcal{A}_{f})$ are complete cotorsion pairs, a cofibrant replacement of $M$ is precisely a special right $\mathcal{A}_{c}$-approximation (special $\mathcal{A}_{c}$-precover) of $M$, and a fibrant replacement of $M$ is a special left $\mathcal{A}_{f}$-approximation (special $\mathcal{A}_{f}$-preenvelope) of $M$.

Suppose $\mathcal{C}$ and $\mathcal{D}$ are model categories, and $(F, G, \varphi): \mathcal{C}\rightarrow \mathcal{D}$ is an adjunction. We call $(F, G, \varphi)$ a Quillen adjunction if $F$ is a left Quillen functor, or equivalently $G$ is a right Quillen functor. That is, $F$ preserves cofibrations and trivial cofibrations, or $G$ preserves fibrations and trivial fibrations. A Quillen adjunction $(F, G, \varphi): \mathcal{C}\rightarrow \mathcal{D}$ is a Quillen equivalence if and only if $(\mathrm{L}F, \mathrm{R}G, \mathrm{R}\varphi): \mathrm{Ho}(\mathcal{C})\rightarrow \mathrm{Ho}(\mathcal{D})$ is an adjoint equivalence of homotopy categories (see e.g. \cite[Proposition 1.3.13]{Hov99}), where the total left derived functor $\mathrm{L}F: \mathrm{Ho}(\mathcal{C})\rightarrow \mathrm{Ho}(\mathcal{D})$ of $F$ is defined to be the composition $\mathrm{Ho}(\mathcal{C})\stackrel{\mathrm{Ho}Q}\longrightarrow \mathrm{Ho}(\mathcal{C}_{c})\stackrel{\mathrm{Ho}F}\longrightarrow \mathrm{Ho}(\mathcal{D})$, and the total right derived functor $\mathrm{R}G: \mathrm{Ho}(\mathcal{D})\rightarrow \mathrm{Ho}(\mathcal{C})$ of $G$ is defined to be the composition $\mathrm{Ho}(\mathcal{D})\stackrel{\mathrm{Ho}R}\longrightarrow \mathrm{Ho}(\mathcal{D}_{f})\stackrel{\mathrm{Ho}G}\longrightarrow  \mathrm{Ho}(\mathcal{C})$.

Our first result shows the equivalence of the stable categories arising from a cotorsion triple. It is easy to see that the identity adjunctions are Quillen adjunction between the model structures $\mathcal{M}^{proj}$ and $\mathcal{M}^{inj}$. We will show that the identity adjunction is indeed a Quillen equivalence, and then, we get the desired equivalence of stable categories $\mathcal{X}/\mathcal{P}\simeq \mathcal{Y}/\mathcal{I}$. We note that an argument for the equivalence $\mathcal{X}/\mathcal{P}\simeq \mathcal{Y}/\mathcal{I}$ is also given in \cite[Chapter VI, Theorem 3.2]{BR07}, by using torsion pairs $(\mathcal{X}/\mathcal{P}, \mathcal{Z}/\mathcal{P})$ and $(\mathcal{Z}/\mathcal{I}, \mathcal{Y}/\mathcal{I})$ in the stable categories $\mathcal{A}/\mathcal{P}$ and $\mathcal{A}/\mathcal{I}$.

\begin{thm}\label{thm 1}
Let $\mathcal{A}$ be a bicomplete abelian category with a complete and hereditary cotorsion triple $(\mathcal{X}, \mathcal{Z}, \mathcal{Y})$. There are equivalences of triangulated categories
$$ \mathcal{X}/\mathcal{P}\simeq \mathrm{Ho}(\mathcal{M}^{proj})\simeq \mathrm{Ho}(\mathcal{M}^{inj})\simeq \mathcal{Y}/\mathcal{I}.$$
Then, the model category $\mathcal{A}$ (with respect to either $\mathcal{M}^{proj}$ or $\mathcal{M}^{inj}$) is stable in the sense of \cite[Definition 7.1.1]{Hov99}.
\end{thm}

\begin{proof}
For a bicomplete abelian category $\mathcal{A}$ with the model structure $\mathcal{M} = [\mathcal{A}_{c}, \mathcal{A}_{tri}, \mathcal{A}_{f}]$, the inclusion functors induce equivalences of categories $\mathrm{Ho}(\mathcal{A}_{cf})\rightarrow \mathrm{Ho}(\mathcal{A}_{c})\rightarrow \mathrm{Ho}(\mathcal{A})$ and $\mathrm{Ho}(\mathcal{A}_{cf})\rightarrow \mathrm{Ho}(\mathcal{A}_{f})\rightarrow \mathrm{Ho}(\mathcal{A})$, where $\mathcal{A}_{cf}=\mathcal{A}_{c}\cap \mathcal{A}_{f}$. Moreover, the inclusion $\mathcal{A}_{cf}\rightarrow \mathcal{A}$ induces an equivalence of categories ${\mathcal{A}_{cf}/\sim} = \mathcal{A}_{cf}/\omega \rightarrow \mathrm{Ho}(\mathcal{A}_{cf})\rightarrow \mathrm{Ho}(\mathcal{A})$, where $f\sim g: M\rightarrow N$ if $g-f$ factors through an object in $\omega = \mathcal{A}_{c}\cap \mathcal{A}_{tri}\cap \mathcal{A}_{f}$. Let $i$ be the inclusion functors. Then there are equivalences of categories $\xymatrix@C=40pt{\mathcal{X}/\mathcal{P}\ar@<0.4ex>[r]^{\mathrm{Ho}i\quad} &\mathrm{Ho}(\mathcal{M}^{proj})\ar@<0.4ex>[l]^{\mathrm{Ho}Q\quad}}$ and  $\xymatrix@C=40pt{\mathcal{Y}/\mathcal{I}\ar@<0.4ex>[r]^{\mathrm{Ho}i\quad} &\mathrm{Ho}(\mathcal{M}^{inj})\ar@<0.4ex>[l]^{\mathrm{Ho}R\quad}}$.
It remains to prove the equivalence $\mathrm{Ho}(\mathcal{M}^{proj})\simeq \mathrm{Ho}(\mathcal{M}^{inj})$.

There is a useful criterion for checking the given Quillen adjunction is a Quillen equivalence, see details in [25, Corollary 1.3.16]. Specifically, we need to show that the identity functor reflects weak equivalences between cofibrant objects in $\mathcal{M}^{proj}=[\mathcal{X}, \mathcal{Z}, \mathcal{A}]$, and for every fibrant object $M$ in $\mathcal{M}^{inj} = [\mathcal{A}, \mathcal{Z}, \mathcal{Y}]$, $QM\rightarrow M$ is a weak equivalence in $\mathcal{M}^{inj}$, where $QM$ is a cofibrant replacement of $M$ with respect to $\mathcal{M}^{proj}$.

Assume that $X_1$ and $X_2$ are cofibrant objects in $\mathcal{M}^{proj}$, i.e. $X_1, X_2\in \mathcal{X}$, and let $f: X_1\rightarrow X_2$ be a weak equivalence in $\mathcal{M}^{inj}$. We need to show $f: X_1\rightarrow X_2$ is also a weak equivalence in $\mathcal{M}^{proj}$, that is, the identity functor reflects weak equivalences between cofibrant objects in $\mathcal{M}^{proj}$. In $\mathcal{M}^{proj}$ we factor $f$ as $f = pi$, where $i: X_1\rightarrow N$ is a trivial cofibration (i.e. a monomorphism with projective cokernel), and $p: N\rightarrow X_2$ is a fibration (i.e. an epimorphism). It follows from the exact sequence
$0\rightarrow X_1\stackrel{i}\rightarrow N\rightarrow \mathrm{Coker}(i)\rightarrow 0$ that $N\in \mathcal{X}$. Let $K = \mathrm{Ker}(p)$. It follows from the exact sequence $0\rightarrow K\rightarrow N\rightarrow X_2\rightarrow 0$ that $K\in \mathcal{X}$. For $K$, there exists an exact sequence $0\rightarrow K\rightarrow Y\rightarrow Z\rightarrow 0$ with $Y\in \mathcal{Y}$ and
$Z\in \mathcal{Z}$. Consider the pushout of $K\rightarrow N$ and $K\rightarrow Y$:
$$\xymatrix@C=20pt@R=20pt{ & 0\ar[d] & 0\ar[d] \\
0 \ar[r]^{}  &K \ar[d] \ar[r] & N \ar[d]_{j}
  \ar[r]^{p} &X_2 \ar@{=}[d]  \ar[r] &0 \\
0 \ar[r] & Y \ar[r] \ar[d] & L\ar[r]^{p'} \ar[d] & X_2 \ar[r] & 0\\
  & Z \ar[d] \ar@{=}[r] & Z\ar[d]\\
  & 0 & 0
  }$$
Note that $i: X_1\rightarrow N$ is also a trivial cofibration in $\mathcal{M}^{inj}$, then $f: X_1\rightarrow X_2$ is a weak equivalence in $\mathcal{M}^{inj}$ if and only if so is  $p: N\rightarrow X_2$.  Moreover, $j:N\rightarrow L$ is a trivial cofibration in $\mathcal{M}^{inj}$, and it yields that $p':L\rightarrow X_2$ is a trivial fibration in $\mathcal{M}^{inj}$, i.e.  $Y\in \mathcal{Y}\cap \mathcal{Z} = \mathcal{I}$. From the left column of the above diagram, we get $K\in \mathcal{Z}$. Hence, $K\in \mathcal{Z}\cap \mathcal{X} = \mathcal{P}$, and it yields that $p: N\rightarrow X_2$ is a trivial fibration,
and then $f=pi$ is a weak equivalence in $\mathcal{M}^{proj}$, as expected.

Let $M$ be a fibrant object in $\mathcal{M}^{inj}$, i.e. $M\in\mathcal{Y}$. In $\mathcal{M}^{proj}=[\mathcal{X}, \mathcal{Z}, \mathcal{A}]$, a cofibrant replacement $g: QM\rightarrow M$ is precisely a special right $\mathcal{X}$-approximation of $M$, i.e. an epic with $K=\mathrm{Ker}(g)\in \mathcal{Z}$. Similarly, there exists an exact sequence $0\rightarrow K\rightarrow Y\rightarrow Z\rightarrow 0$ with $Y\in \mathcal{Y}$ and
$Z\in \mathcal{Z}$. Consider the pushout of $K\rightarrow QM$ and $K\rightarrow Y$:
$$\xymatrix@C=20pt@R=20pt{ & 0\ar[d] & 0\ar[d] \\
0 \ar[r]^{}  &K \ar[d] \ar[r] & QM \ar[d]_{i}
  \ar[r]^{g} &M \ar@{=}[d]  \ar[r] &0 \\
0 \ar[r] & Y \ar[r] \ar[d] & P\ar[r]^{p} \ar[d] & M \ar[r] & 0\\
  & Z \ar[d] \ar@{=}[r] & Z\ar[d]\\
  & 0 & 0
  }$$
Note that $Y\in \mathcal{Y}\cap\mathcal{Z} = \mathcal{I}$, and then $p$ is a trivial fibration in $\mathcal{M}^{inj}$; it is obvious that $i$ is a trivial cofibration in $\mathcal{M}^{inj}$. Thus,  $g = pi$ is a weak equivalence in $\mathcal{M}^{inj}$.

Then the identity adjunction $\xymatrix@C=20pt{\mathcal{M}^{proj}\ar@<0.4ex>[r]^{} &\mathcal{M}^{inj}\ar@<0.4ex>[l]^{}}$ is a Quillen equivalence. This implies the desired adjoint equivalence of homotopy categories $\xymatrix@C=40pt{\mathrm{Ho}(\mathcal{M}^{proj})\ar@<0.4ex>[r]^{\mathrm{L}\mathrm{id}} &\mathrm{Ho}(\mathcal{M}^{inj})\ar@<0.4ex>[l]^{\mathrm{R}\mathrm{id}}}$. The homotopy categories and the stable categories must be triangulated categories as shown in \cite[Chapter 7]{Hov99}. This completes the proof.
\end{proof}

The morphisms in the homotopy categories are characterized as following.

\begin{prop}\label{prop 6}
Let $\mathcal{A}$ be a bicomplete abelian category with a complete and hereditary cotorsion triple $(\mathcal{X}, \mathcal{Z}, \mathcal{Y})$,
and $M$, $N$ be objects in $\mathcal{A}$. Then $\mathrm{Hom}_{\mathrm{Ho}(\mathcal{A})}(M, N)\simeq \mathrm{Hom}_{\mathcal{A}}(M, Y^N)/\sim$,
where $Y^N$ is a special left $\mathcal{Y}$-approximation of $N$ and maps $f\sim g: M\rightarrow Y^{N}$ if $g-f$ factors through an injective object. Dually, $\mathrm{Hom}_{\mathrm{Ho}(\mathcal{A})}(M, N)\simeq \mathrm{Hom}_{\mathcal{A}}(X_M, N)/\sim$, where
$X_M$ is a special right $\mathcal{X}$-approximation of $M$ and maps $f\sim g: X_{M}\rightarrow N$ if $g-f$ factors through a projective object.
\end{prop}

\begin{proof}
In any model category $\mathcal{A}$, it follows from \cite[Theorem 1.2.10]{Hov99} that
$\mathrm{Hom}_{\mathrm{Ho}(\mathcal{A})}(M, N)\simeq \mathrm{Hom}_{\mathcal{A}}(QM, RN)/\sim$,
where $\sim$ are the left homotopy relations (equivalently, right homotopy relations), and are equivalence relations on
$\mathrm{Hom}_{\mathcal{A}}(QM, RN)$. See \cite{Hov99} and \cite{DS95} for details of cylinder object, path object, left homotopy and
right homotopy.

Consider the injective model structure $[\mathcal{A}, \mathcal{Z}, \mathcal{Y}]$ on $\mathcal{A}$.  For $N$, $i: N\rightarrow Y^N$ is a trivial cofibration if and only if it is a special left $\mathcal{Y}$-approximation, and then $RN = Y^N$. Note that $QM = M$ since $M$ is cofibrant itself. It remains to determine the homotopy relation.

Let $j: M\rightarrow I$ be a monomorphism with $I$ injective. Then
$$M\oplus M \stackrel{\alpha}\longrightarrow M\oplus I\oplus I \stackrel{\beta}\longrightarrow M$$
is a factorization of the folding map, where $\alpha = \left(
                             \begin{array}{cc}
                             1_{M} & 1_{M} \\
                             j & 0 \\
                             0 & j \\
                             \end{array}
                            \right)$
is a monomorphism (i.e. a cofibration in the injective model structure),
$\beta = (1_{M}, 0, 0)$ is an epimorphism with injective kernel (i.e. a trivial fibration in the injective model structure).
So $M\oplus I\oplus I$ is a cylinder object for $M$. Then
$f, g: M\rightarrow Y^N$ are left homotopic if and only if there is a map $H = (r, s, t)$ such that $H\alpha = (f, g)$, i.e. $f = r + sj$
and $g = r + tj$. If there is such a map, then $g-f = (t-s)j$, so $g-f$ factors through the injective object $I$.

Conversely, assume $g-f$ factors through some injective object $I'$. Then the map $M\rightarrow I'$ extends to a map $I\rightarrow I'$, so
in fact $g-f = tj$ for some $t:I\rightarrow Y$. Then $f = f + 0j$ and $g= f+tj$. Thus $f$ and $g$ are left homotopic.

The second statement can be proved dually by considering the projective model structure on $\mathcal{A}$, and constructing the path
object and right homotopy.
\end{proof}

It is useful to know when two objects become equivalent in the stable categories.

\begin{prop}\label{prop 7}
Let $\mathcal{A}$ be a bicomplete abelian category with a complete and hereditary cotorsion triple $(\mathcal{X}, \mathcal{Z}, \mathcal{Y})$,
and $M$, $N$ be objects in $\mathcal{A}$. Then $M$ and $N$ are equivalent in the stable category $\mathcal{X}/\mathcal{P}$
if and only if
$X_M\oplus P \cong X_N\oplus F$, where $X_M$ and $X_N$ are special right $\mathcal{X}$-approximations of $M$ and $N$ respectively,
$P$ and $F$ are projective objects.
Dually,  $M$ and $N$ are equivalent in the stable category $\mathcal{Y}/\mathcal{I}$ if and only if
$Y^M\oplus I \cong Y^N\oplus J$, where $Y^M$ and $Y^N$ are special left $\mathcal{Y}$-approximations of
$M$ and $N$ respectively, $I$ and $J$ are injective objects.
\end{prop}

\begin{proof}
Note that $M$ and $N$ are equivalent in the category $\mathcal{X}/\mathcal{P}$ if and only if the morphism $f: QM\rightarrow QN$ is a weak equivalence
with respect to the projective model structure. Factor $f$ into a trivial cofibration $f':QM\rightarrow L$
followed by a trivial fibration $f'': L\rightarrow QN$. Then $f'$ is a split monomorphism with $P = \mathrm{Coker}(f')$ projective, and
$f''$ is a split epimorphism since $F = \mathrm{Ker}(f'')\in \mathcal{Z}$ and $QN\in \mathcal{X}$.
So $QM\oplus P\cong L\simeq QN \oplus F$. Moreover, since $F$ is in $\mathcal{X}$ as well, then $F$ is projective.
Note that the cofibrant replacement $QM$ for any object $M\in\mathcal{A}$
is precisely a special right $\mathcal{X}$-approximation $X_M$ of $M$, the desired result follows.

Similarly, we can prove the equivalence in the stable category $\mathcal{Y}/\mathcal{I}$
by considering the injective model structure on $\mathcal{A}$.
\end{proof}

\section {\bf Homological dimensions with respect to cotorsion triples}

We start with the following definition.

\begin{df}\label{def 2}
Given a subcategory $\mathcal{C}$ of an abelian category  $\mathcal{A}$.
For any object $M\in \mathcal{A}$, the $\mathcal{C}$-projective dimension $\mathcal{C}$-$\mathrm{pd}(M)$ of $M$ is defined to be
the smallest integer $n \geq 0$ such that $\mathrm{Ext}_{\mathcal{A}}^{n+i}(M, C) = 0$ for any
$C\in \mathcal{C}$ and any $i \geq 1$. Dually, the $\mathcal{C}$-injective dimension $\mathcal{C}$-$\mathrm{id}(M)$ of $M$ is defined to be
the smallest integer $n \geq 0$ such that $\mathrm{Ext}_{\mathcal{A}}^{n+i}(C, M) = 0$ for any
$C\in \mathcal{C}$ and any $i \geq 1$. If there is no such $n$ existed, let $\mathcal{C}$-$\mathrm{pd}(M)$ ($\mathcal{C}$-$\mathrm{id}(M)$)
be $\infty$.

The global $\mathcal{C}$-projective (resp. $\mathcal{C}$-injective) dimension  $\mathcal{C}$-$\mathrm{g.pd}(\mathcal{A})$
(resp. $\mathcal{C}$-$\mathrm{g.id}(\mathcal{A})$) is defined to be the supremum of the  $\mathcal{C}$-projective (resp. $\mathcal{C}$-injective)
dimension for all objects in $\mathcal{A}$.
\end{df}

Recall that for a contravariant finite (or, precoving) subcategory $\mathcal{X}\subseteq \mathcal{A}$ and an object $M\in \mathcal{A}$,
an $\mathcal{X}$-resolution of $M$ is a complex
$\cdots \rightarrow X_1\stackrel{d_1}\rightarrow X_0\stackrel{d_0}\rightarrow M\rightarrow 0$ with each $X_{i}\in \mathcal{X}$
such that it is exact after applying $\mathrm{Hom}_{\mathcal{A}}(X, -)$ for each $X\in\mathcal{X}$; this is
equivalent to that each induced homomorphism $X_{i}\rightarrow \mathrm{Ker}(d_{i-1})$ for any $i \geq 0$ (set $M = \mathrm{Ker}(d_{-1})$)
is a right $\mathcal{X}$-approximation. We denote sometimes the $\mathcal{X}$-resolution by $\mathbf{X}^{\bullet}\rightarrow M$, where
$\mathbf{X}^{\bullet} = \cdots \rightarrow X_1\stackrel{d_1}\rightarrow X_0 \rightarrow 0$ is the deleted $\mathcal{X}$-resolution of $M$.
Dually, for a covariant finite (or, preenveloping) subcategory $\mathcal{Y}$, one has the notion of $\mathcal{Y}$-coresolution and then the notion of deleted $\mathcal{Y}$-coresolution.

\begin{prop}\label{prop 8}
Let $(\mathcal{X}, \mathcal{Z}, \mathcal{Y})$ be a complete and hereditary cotorsion triple in $\mathcal{A}$.
For any object $M\in\mathcal{A}$, we have  $\mathcal{X}$-$\mathrm{id}(M) = \mathcal{Y}$-$\mathrm{pd}(M) = 0$
or $\infty$. Consequently, $\mathcal{X}$-$\mathrm{g.id}(\mathcal{A}) = \mathcal{Y}$-$\mathrm{g.pd}(\mathcal{A}) = 0$
 or $\infty$.
\end{prop}

\begin{proof}
Assume that $\mathcal{X}$-$\mathrm{id}(M)$ is finite, say $n$. For $M$, there is a  $\mathcal{Z}$-coresolution
$0\rightarrow M\rightarrow Z^{0}\rightarrow \cdots\rightarrow Z^{n-1}\rightarrow Z^{n}\rightarrow\cdots$.
Set $M = Z^{-1}$ and $L^{i} = \mathrm{Coker}(Z^{i-2}\rightarrow Z^{i-1})$ for all $i \geq 1$. By dimension shifting, we
have $$0 = \mathrm{Ext}_{\mathcal{A}}^{n+1}(X, M)\simeq \mathrm{Ext}_{\mathcal{A}}^{n}(X, L^{1})
\simeq \cdots \simeq \mathrm{Ext}_{\mathcal{A}}^{1}(X, L^{n})$$ for any $X\in \mathcal{X}$. This implies $L^{n}\in \mathcal{Z}$.
Note that $\mathcal{Z}$ is thick, we have inductively $L^{n}$, $L^{n-1}$, $\cdots$, and $M$ are in $\mathcal{Z}$. Then
$\mathrm{Ext}_{\mathcal{A}}^{1}(X, M) = 0$ for any $X\in \mathcal{X}$ and hence
$\mathcal{X}$-$\mathrm{id}(M) = 0$. Similarly, if $\mathcal{Y}$-$\mathrm{pd}(M)$ is finite, then $\mathcal{Y}$-$\mathrm{pd}(M) = 0$.
The desired assertion follows.
\end{proof}

Let $(\mathcal{X}, \mathcal{Z}, \mathcal{Y})$ be a complete and hereditary cotorsion triple in $\mathcal{A}$, and $M$ be any
object of $\mathcal{A}$. It follows immediately from Propositions \ref{prop 3} and \ref{prop 4} that:
each $\mathcal{X}$-resolution of $M$ remains exact after applying  $\mathrm{Hom}_{\mathcal{A}}(-, Y)$
for each $Y\in\mathcal{Y}$, and similarly, each $\mathcal{Y}$-coresolution of $M$ remains exact after applying
$\mathrm{Hom}_{\mathcal{A}}(X, -)$ for each $X\in\mathcal{X}$. Then $\mathrm{Hom}_{\mathcal{A}}(-, -)$ is right balanced
by $\mathcal{X}\times \mathcal{Y}$ in the sense of \cite[Definition 8.2.13]{EJ00}.
Let us remark that the pair $(\mathcal{X}, \mathcal{Y})$ is a balanced pair
in the sense of \cite[Definition 1.1]{Chen10}.

An advantage of such balanced properties is that the cohomology groups computed by $\mathcal{X}$-resolution of the first variable
and by $\mathcal{Y}$-coresolution of the second variable are isomorphic.
By \cite[Theorem 8.2.14]{EJ00} or \cite[Lemma 2.1]{Chen10}, we have the following definition.

\begin{df}\label{def 3}
Let $(\mathcal{X}, \mathcal{Z}, \mathcal{Y})$ be a complete and hereditary cotorsion triple in $\mathcal{A}$.
For any objects $M, N$ of $\mathcal{A}$ and any $i \geq 0$, we define
$$\mathrm{Ext}^{i}_{\mathcal{XY}}(M, N) := H^i\mathrm{Hom}_{\mathcal{A}}(\mathbf{X}^{\bullet}, N)
\cong H^i\mathrm{Hom}_{\mathcal{A}}(M, \mathbf{Y}^{\bullet}),$$
where $\mathbf{X}^{\bullet}\rightarrow M$ is an $\mathcal{X}$-resolution of $M$,
and $N\rightarrow \mathbf{Y}^{\bullet}$  is a $\mathcal{Y}$-coresolution of $N$.
\end{df}

By a version of comparison theorem, both $\mathcal{X}$-resolution and  $\mathcal{Y}$-coresolution are unique up to homotopy, and then
the functors $\mathrm{Ext}^{i}_{\mathcal{XY}}(-, -)$ are well defined. We should call such functors relative extension functors
(with respect to the cotorsion triple $(\mathcal{X}, \mathcal{Z}, \mathcal{Y})$).

It is easy to see that for any objects $M, N\in \mathcal{A}$,
$\mathrm{Ext}^{0}_{\mathcal{XY}}(M, N)\simeq \mathrm{Hom}_{\mathcal{A}}(M, N)\simeq\mathrm{Ext}^{0}_{\mathcal{A}}(M, N)$.
We want to compare $\mathrm{Ext}^{i}_{\mathcal{XY}}(-, -)$ with $\mathrm{Ext}^{i}_{\mathcal{A}}(-, -)$ for $i > 0$.

\begin{prop}\label{prop 9}
Let $(\mathcal{X}, \mathcal{Z}, \mathcal{Y})$ be a complete and hereditary cotorsion triple in $\mathcal{A}$,
and $M, N$ be objects in $\mathcal{A}$. If either $M$ or $N$ is in $\mathcal{Z}$, then  for any $i> 0$ we have
$$\mathrm{Ext}^{i}_{\mathcal{XY}}(M, N)\cong \mathrm{Ext}^{i}_{\mathcal{A}}(M, N).$$
\end{prop}

\begin{proof}
Assume that $N\in \mathcal{Z}$. There is a $\mathcal{Y}$-coresolution $0\rightarrow N\rightarrow Y^{0}\rightarrow \cdots
\rightarrow Y^{n-1}\rightarrow Y^{n}\rightarrow \cdots$. Note that $\mathrm{Coker}(N\rightarrow Y^{0})\in \mathcal{Z}$.
Then $Y^{0}\in \mathcal{Z} \cap \mathcal{Y}$. It follows from Proposition \ref{prop 2} that $Y^{0}$ is injective.
Inductively, we have each $Y^{j}$ is injective for $j\geq 0$. Then $N\rightarrow \mathbf{Y}^{\bullet}$ is also an injective coresolution.
Hence, $\mathrm{Ext}^{i}_{\mathcal{XY}}(M, N)\cong \mathrm{Ext}^{i}_{\mathcal{A}}(M, N)$.
Similarly, we can prove $\mathrm{Ext}^{i}_{\mathcal{XY}}(M, N)\cong \mathrm{Ext}^{i}_{\mathcal{A}}(M, N)$ when $M\in \mathcal{Z}$.
\end{proof}

We have the following characterization on homological dimensions with respect to a cotorsion triple by a standard argument in classical homological algebra. So we omit the proof.

\begin{thm}\label{thm 2}
Let $(\mathcal{X}, \mathcal{Z}, \mathcal{Y})$ be a complete and hereditary cotorsion triple in $\mathcal{A}$.
For any object $M\in \mathcal{A}$ and nonnegative integer $n$, the following are equivalent:\\
\indent $\mathrm{(1)}$  $\mathcal{Z}\text{-}\mathrm{pd}(M) \leq n$.\\
\indent $\mathrm{(2)}$  $\mathrm{Ext}_{\mathcal{A}}^{n+i}(M, Z) = 0$ for any $Z\in \mathcal{Z}$ and $i \geq 1$.\\
\indent $\mathrm{(3)}$  $\mathrm{Ext}_{\mathcal{XY}}^{n+i}(M, N) = 0$ for any $N\in \mathcal{A}$  and $i \geq 1$.\\
\indent $\mathrm{(4)}$  For any exact sequence
$\cdots \rightarrow X_{n}\rightarrow X_{n-1}\rightarrow \cdots \rightarrow X_0\rightarrow M\rightarrow 0$ with each $X_{i}$ in $\mathcal{X}$,
$K_{n} = \mathrm{Ker}(X_{n-1}\rightarrow X_{n-2})\in \mathcal{X}$.
\end{thm}

Dually, we have the following.

\begin{thm}\label{thm 3}
Let $(\mathcal{X}, \mathcal{Z}, \mathcal{Y})$ be a complete and hereditary cotorsion triple in $\mathcal{A}$.
For any object $N\in \mathcal{A}$ and nonnegative integer $n$, the following are equivalent:\\
\indent $\mathrm{(1)}$  $\mathcal{Z}\text{-}\mathrm{id}(N) \leq n$.\\
\indent $\mathrm{(2)}$  $\mathrm{Ext}_{\mathcal{A}}^{n+i}(Z, N) = 0$ for any $Z\in \mathcal{Z}$ and $i \geq 1$.\\
\indent $\mathrm{(3)}$  $\mathrm{Ext}_{\mathcal{XY}}^{n+i}(M, N) = 0$ for any $M\in \mathcal{A}$  and $i \geq 1$.\\
\indent $\mathrm{(4)}$  For any exact sequence
$0 \rightarrow N\rightarrow Y^{0}\rightarrow \cdots \rightarrow  Y^{n-1}\rightarrow Y^n\rightarrow\cdots$ with each $Y^{i}$ in $\mathcal{Y}$,
$L^{n} = \mathrm{Coker}(Y^{n-1}\rightarrow Y^{n})\in \mathcal{Y}$.
\end{thm}

Recall that the $\mathcal{X}$-resolution dimension $\mathcal{X}$-res.dim$(M)$ of an object $M$ \cite{Chen10}
is defined to be the minimal length of
$\mathcal{X}$-resolutions of $M$. By comparing Theorem \ref{thm 2} with \cite[Lemma 2.4]{Chen10}, we have
$\mathcal{Z}$-pd$(M) = \mathcal{X}$-res.dim$(M)$. Dually, $\mathcal{Z}$-injective dimension of $M$
equals to $\mathcal{Y}$-coresolution dimension of $M$.
It is not difficult to get the following, which implies ``the balance of $\mathcal{X}$ and $\mathcal{Y}$''
in the sense of homological dimensions.

\begin{cor}\label{cor 1}
Let $(\mathcal{X}, \mathcal{Z}, \mathcal{Y})$ be a complete and hereditary cotorsion triple in $\mathcal{A}$.
Then we have $\mathcal{Z}\text{-}\mathrm{g.pd}(\mathcal{A}) = \mathcal{Z}\text{-}\mathrm{g.id}(\mathcal{A})$
(or, $\mathcal{X}\text{-}\mathrm{res.dim}(\mathcal{A}) = \mathcal{Y}\text{-}\mathrm{cores.dim}(\mathcal{A}))$.
\end{cor}

Moreover, we have the following.

\begin{cor}\label{cor 2}
Let $(\mathcal{X}, \mathcal{Z}, \mathcal{Y})$ be a complete and hereditary cotorsion triple in $\mathcal{A}$.
Then the following are equivalent for a nonnegative integer $n$:\\
\indent$\mathrm{(1)}$ $\mathcal{Z}\text{-}\mathrm{g.pd}(\mathcal{A}) = \mathcal{Z}\text{-}\mathrm{g.id}(\mathcal{A}) \leq n$.\\
\indent$\mathrm{(2)}$  Each $Z\in \mathcal{Z}$ has projective dimension at most $n$.\\
\indent$\mathrm{(3)}$  Each $Z\in \mathcal{Z}$ has injective dimension at most $n$.
\end{cor}

\section {\bf Cotorsion triples of complexes and homotopy equivalences}

Let $(\mathcal{X}, \mathcal{Y})$ be a balanced pair in the sense of \cite[Definition 1.1]{Chen10}, which arises naturally from a cotorsion pair. Chen \cite{Chen10} proved that if the global $\mathcal{X}$-resolution dimension (equally, global $\mathcal{Y}$-coresolution dimension) is finite, then the relative derived categories $\mathrm{D}_{\mathcal{X}}(R)$ and $\mathrm{D}_{\mathcal{Y}}(R)$ are equivalent, and then from natural composition equivalent functors $\mathrm{K}(\mathcal{X})\stackrel{inc}\rightarrow \mathrm{K}(R)\rightarrow \mathrm{D}_{\mathcal{X}}(R)$ and $\mathrm{K}(\mathcal{Y})\stackrel{inc}\rightarrow \mathrm{K}(R)\rightarrow \mathrm{D}_{\mathcal{Y}}(R)$, there is a homotopy equivalence
$\mathrm{K}(\mathcal{X})\simeq \mathrm{K}(\mathcal{Y})$; see \cite[Theorem A]{Chen10}.

In this section, we firstly lift a cotorsion triple in an abelian category to associated cotorsion triples of chain complexes. As an application, we can prove the above equivalence of homotopy categories by Quillen equivalence.

The following notations were introduced in \cite{Gil04} and \cite{Gil08}.

\begin{df} \label{df 4}
Given a class of objects $\mathcal{X}$ in an abelian category $\mathcal{A}$, we define the following classes of chain complexes in $\mathrm{Ch}(\mathcal{A})$.\\
\indent $(1)$ $dw\widetilde{\mathcal{X}}$ is the class of all chain complexes $X$ with $X_n\in\mathcal{X}$.\\
\indent $(2)$ $\widetilde{\mathcal{X}}$ is the class of all exact chain complexes $X$ with cycles $Z_{n} X\in\mathcal{X}$.\\
\indent $(3)$ If $(\mathcal{X}, \mathcal{Y})$ is a cotorsion pair in $\mathcal{A}$, then we will denote $\widetilde{\mathcal{X}}^{\bot}$ by $dg\widetilde{\mathcal{Y}}$ and $^{\bot}\widetilde{\mathcal{Y}}$ by $dg\widetilde{\mathcal{X}}$.
\end{df}

The suspension of the complex $X$, denoted $\Sigma X$, is the complex given by  $(\Sigma X)_{n} = X_{n-1}$ and $d_{\Sigma X}^{n} = - d_{n-1}$. The complex $\Sigma (\Sigma X)$ is denoted $\Sigma^{2}X$, and inductively, we define
$\Sigma^{n}X$ for all $n\in \mathbb{Z}$.
Given two complexes $X$ and $Y$, let $\mathcal{H}om(X, Y)$ denote the complex of $\mathbb{Z}$-modules with $n$th component
$\mathcal{H}om(X, Y)_{n} = \prod_{k\in \mathbb{Z}}\textrm{Hom}(X_{k}, Y_{k+n})$ and differential
$(\delta_{n}f)_{k} = d_{k+n}^{Y}f_{k} - (-1)^{n}f_{k-1}d_{k}^{X}$
for morphisms $f_{k}: X_{k}\rightarrow Y_{k+n}$.

Let $\mathcal{A}$ be an abelian category with a complete and hereditary cotorsion triple $(\mathcal{X}, \mathcal{Z}, \mathcal{Y})$. By Corollary \ref{cor 1}, we say $\mathcal{A}$ is of finite global $\mathcal{Z}$-dimension if $\mathcal{Z}\text{-}\mathrm{g.pd}(\mathcal{A}) = \mathcal{Z}\text{-}\mathrm{g.id}(\mathcal{A}) < \infty$ (i.e.
$\mathcal{X}\text{-res.dim}(\mathcal{A}) = \mathcal{Y}\text{-cores.dim}(\mathcal{A})< \infty$). The following generalizes \cite[Proposition 7.4, 7.8]{Gil16}.

\begin{lem}\label{lem 2}
Let $\mathcal{A}$ be an abelian category with a complete and hereditary cotorsion triple $(\mathcal{X}, \mathcal{Z}, \mathcal{Y})$. If $\mathcal{A}$ is of finite global $\mathcal{Z}$-dimension, then there is a complete and hereditary cotorsion triple $(dw\widetilde{\mathcal{X}}, (dw\widetilde{\mathcal{X}})^{\perp} = {^{\perp}(dw\widetilde{\mathcal{Y}})},  dw\widetilde{\mathcal{Y}})$ in $\mathrm{Ch}(\mathcal{A})$.
\end{lem}

\begin{proof}
By \cite{Gil04}, every complete and hereditary cotorsion pair $(\mathcal{X}, \mathcal{Z})$ in $\mathcal{A}$ induces two complete and hereditary cotorsion pairs $(dg\widetilde{\mathcal{X}}, \widetilde{\mathcal{Z}})$ and $(\widetilde{\mathcal{X}}, dg\widetilde{\mathcal{Z}})$ in the category of chain complexes $\mathrm{Ch}(\mathcal{A})$. Then there is a complete and hereditary cotorsion triple $(dg\widetilde{\mathcal{X}}, \widetilde{\mathcal{Z}}, dg\widetilde{\mathcal{Y}})$ in $\mathrm{Ch}(\mathcal{A})$ induced by $(\mathcal{X}, \mathcal{Z}, \mathcal{Y})$. We claim that $dw\widetilde{\mathcal{X}} = dg\widetilde{\mathcal{X}}$ and  $dw\widetilde{\mathcal{Y}} = dg\widetilde{\mathcal{Y}}$.

It follows from \cite{Gil04} that $dg\widetilde{\mathcal{X}} \subseteq dw\widetilde{\mathcal{X}}$. We prove now the inverse inclusion. Let $X$ be a complex in $dw\widetilde{\mathcal{X}}$. There exists a quasi-isomorphism $f: P\rightarrow X$
with $P$ dg-projective. We denote by $M = \mathrm{Cone}(f)$ the mapping cone of $f$. Then $M$ is an exact complex with items $M_i=P_{i-1}\oplus X_i$ being in $\mathcal{X}$. For any cycle $Z_{i}M$, there exist exact sequences
$$0\longrightarrow Z_{i}M\longrightarrow M_{i}\longrightarrow M_{i-1}\longrightarrow \cdots\longrightarrow M_{i-k}\longrightarrow Z_{i-k-1}M\longrightarrow 0.$$
Since $\mathcal{Z}\text{-g.pd}(\mathcal{A})$ is finite, we have $Z_{i}M\in \mathcal{X}$ by Theorem \ref{thm 2}. This implies that
$M \in \widetilde{\mathcal{X}}$, and then by \cite[Theorem 3.12]{Gil04} $M\in dg\widetilde{\mathcal{X}}$.

Now consider the short exact sequence $0\rightarrow X\rightarrow M\rightarrow \Sigma P\rightarrow 0$. For any complex $Z\in \widetilde{\mathcal{Z}}$, we get an exact sequence $0\rightarrow \mathcal{H}om(\Sigma P, Z)\rightarrow \mathcal{H}om(M, Z)\rightarrow \mathcal{H}om(X, Z)\rightarrow 0$, where the first two items are exact complexes. Then the complex $\mathcal{H}om(X, Z)$ is exact. Hence, it follows from \cite[Proposition 3.6]{Gil04} that $X\in dg\widetilde{\mathcal{X}}$. Similarly, we can prove $dw\widetilde{\mathcal{Y}} = dg\widetilde{\mathcal{Y}}$. This completes the proof.
\end{proof}

We use $\widehat{\mathcal{Z}}$ to denote $(dw\widetilde{\mathcal{X}})^{\perp} = {^{\perp}(dw\widetilde{\mathcal{Y}})}$. The following is immediate.

\begin{cor}\label{cor 3}
Let $\mathcal{A}$ be an abelian category with a complete and hereditary cotorsion triple $(\mathcal{X}, \mathcal{Z}, \mathcal{Y})$. If $\mathcal{A}$ is of finite global $\mathcal{Z}$-dimension, then on the category $\mathrm{Ch}(\mathcal{A})$, there is a projective model structure $\widetilde{\mathcal{M}}^{proj} = [dw\widetilde{\mathcal{X}}, \widehat{\mathcal{Z}}, \mathrm{Ch}(\mathcal{A})]$, and an injective model structure $\widetilde{\mathcal{M}}^{inj} = [\mathrm{Ch}(\mathcal{A}), \widehat{\mathcal{Z}}, dw\widetilde{\mathcal{Y}}]$.
\end{cor}

Let $f, g : M\rightarrow N$ be two chain maps. By Proposition \ref{prop 6}, if $M$ is in $dw\widetilde{\mathcal{X}}$ and $N$ is any complex, then $f$ and $g$ are formally homotopic in the above projective model structure, if and only if $g-f$ factors through a projective complex. Dually, if $M$ is any complex and $N$ is in $dw\widetilde{\mathcal{Y}}$, then $f$ and $g$ are formally homotopic in the above injective model structure, if and only if $g-f$ factors through an injective complex. But both projective and injective complexes are contractible (split, or homotopically trivial in some literature), and it follows that in these situations the two maps are homotopic if and only if they are chain homotopic in the usual sense.

We now in a position to state the result on equivalences of homotopy categories, which follows by Theorem \ref{thm 1} and Lemma \ref{lem 2}.  We remark that the equivalence $\mathrm{K}(\mathcal{X})\simeq \mathrm{K}(\mathcal{Y})$ also appears in \cite[Theorem A]{Chen10}, where $(\mathcal{X}, \mathcal{Y})$ are supposed to be a balanced pair in the sense of \cite[Definition 1.1]{Chen10}, which arises naturally from a cotorsion triple.

\begin{thm}\label{thm 4}
Let $\mathcal{A}$ be a bicomplete abelian category with a complete and hereditary cotorsion triple $(\mathcal{X}, \mathcal{Z}, \mathcal{Y})$. Assume that $\mathcal{A}$ is of finite global $\mathcal{Z}$-dimension. Then there exist equivalences of categories
$\mathrm{Ho}(\widetilde{\mathcal{M}}^{proj})\simeq\mathrm{K}(\mathcal{X})\simeq \mathrm{K}(\mathcal{Y})\simeq \mathrm{Ho}(\widetilde{\mathcal{M}}^{inj})$.
\end{thm}

\section {\bf Cotorsion triples of Gorenstein modules}

Recall that a ring $R$ is Iwanaga-Gorenstein if it is two-sided noetherian and the regular module $R$ has finite
injective dimension on both sides. Following \cite{Bel00} a ring $R$ is left-Gorenstein provided that any left $R$-module
has finite projective dimension if and only if it has finite injective dimension. By \cite[Definition 3.9]{ZAD14}, a commutative noetherian ring $R$ of finite Krull dimension such that $\mathcal{GP}^{\perp} = {^{\perp}}\mathcal{GI}$ is called virtually Gorenstein.

\begin{lem}\label{lem 3}
There is a complete and hereditary cotorsion triple $(\mathcal{GP}, \mathcal{W}, \mathcal{GI})$ in the category $R$-Mod, if one of the following conditions holds:\\
\indent $(1)$ $R$ is Iwanaga-Gorenstein.\\
\indent $(2)$ $R$ is left-Gorenstein.\\
\indent $(3)$ $R$ is virtually Gorenstein.
\end{lem}

The first assertion is true by \cite{EJ00} or \cite{Hov02}, and the second one comes from \cite{Bel00} or \cite{Chen10}. In the first two cases, $\mathcal{W}$ is the class of modules with finite projective dimension. When $R$ is virtually Gorenstein, it follows from \cite[Theorem A.1]{MS11} or \cite[Section 5.1]{Gil16} that $(\mathcal{GP}, \mathcal{GP}^{\perp})$ is a complete cotorsion pair, and one gets a complete cotorsion pair $(^{\perp}\mathcal{GI}, \mathcal{GI})$ by \cite[Theorem 7.12]{Kra05}. Hence (3) holds. Moreover, Hovey \cite{Hov02} proved that over an Iwanaga-Gorenstein ring the cotorsion pair $(\mathcal{GP}, \mathcal{W})$ is cogenerated by a set, and recently we generalized this result to left-Gorenstein rings by a different method \cite{Ren18}. This has a particular interest since in this case the corresponding projective model structure is cofibrantly generated.

It is worth noting that: (i) Iwanaga-Gorenstein rings are left-Gorenstein (by \cite[Corollary 6.11]{Bel00} or \cite[Chapter 9]{EJ00}), while the converse is not true in general (see \cite{EEGRI}); (ii) every finite Krull dimensional Iwanaga-Gorenstein ring is virtually Gorenstein, and there are virtually Gorenstein rings, for example any local ring $(R, \mathfrak{m})$ with $\mathfrak{m}^{2} = 0$, which are not Iwanaga-Gorenstein; see \cite[Example 3.13]{ZAD14} or \cite[Proposition 6.1 and Remark 6.5]{IK06}.

The Gorenstein cotorsion triple over an Iwanaga-Gorenstein ring can be generalized to the coherent analogy. Gillespie named a module Ding projective if it equals to the cycle $Z_{0}P$ for an exact complex $P$ of projective modules which remains exact after applying $\mathrm{Hom}_{R}(-, F)$ for any flat module $F$, and there is a dual notion of Ding injective modules \cite[Definition 3.2, 3.6]{Gil10}. $R$ is called a Ding-Chen ring it it is two-sided coherent and the regular module $R$ has finite FP-injective dimension on both sides. By  \cite[Corollary 4.5, 4.6]{Gil10}, we have the following.

\begin{lem}\label{lem 4}
If $R$ is a Ding-Chen ring, then there is a complete and hereditary cotorsion triple $(\mathcal{DP}, \mathcal{W}, \mathcal{DI})$ in $R$-Mod, where $\mathcal{DP}$ and $\mathcal{DI}$ are classes of Ding projective and Ding injective modules respectively, $\mathcal{W}$ is the class of modules with finite flat dimension.

Especially, if $R$ is an IF-ring, then there is a complete and hereditary cotorsion triple $(^{\perp}\mathcal{F}, \mathcal{F}, \mathcal{F}^{\perp})$
induced by the subcategory $\mathcal{F}$ of all flat $R$-modules.
\end{lem}

It is well-known that the stable category $\mathcal{GP}/\mathcal{P}$ is closely related to the singularity category. The following is direct from Theorem \ref{thm 1}.

\begin{cor}\label{cor 4}
If $R$ is a left-Gorenstein or virtually Gorenstein ring, then there is an equivalence of categories $\mathcal{GP}/\mathcal{P}\simeq\mathcal{GI}/\mathcal{I}$. If $R$ is a Ding-Chen ring, then $\mathcal{DP}/\mathcal{P}\simeq\mathcal{DI}/\mathcal{I}$.
\end{cor}

For any ring $R$, Bennis and Mahdou proved an equality (\cite[Theorem 1.1]{BM10}):
\begin{center}$\mathrm{sup}\{\mathrm{G}\text{-}\mathrm{proj.dim}_{R}(M) \mid M\in R\text{-}\mathrm{Mod} \} = \mathrm{sup}\{\mathrm{G}\text{-}\mathrm{inj.dim}_{R}(M) \mid M\in R\text{-}\mathrm{Mod} \},$\end{center}
and they named the common value of this equality, denoted by $\textrm{G-gldim}(R)$, the Gorenstein global dimension of $R$.
When $R$ is left-Gorenstein, we consider the cotorsion triple $(\mathcal{GP}, \mathcal{W}, \mathcal{GI})$, and it follows from Proposition \ref{prop 3} and \ref{prop 4} that Gorenstein projective modules are ``projective'' with respect to $\mathcal{GI}$, and Gorenstein injective modules are ``injective'' with respect to $\mathcal{GP}$. Moreover, it follows from Theorem \ref{thm 2} and \ref{thm 3} that we can characterize Gorenstein projective and Gorenstein injective dimensions (i.e. $\mathcal{W}$-projective and $\mathcal{W}$-injective dimensions) of modules by vanishing of relative derived functors; this generalizes some results in \cite[Theorem 4.2, 2(a)]{AM02} and \cite[Proposition 11.5.7]{EJ00}. Note that $R$ is left-Gorenstein if $R$ has finite Gorenstein global dimension. Therefore we can give a new proof for Bennis and Mahdou's equality by vanishing of relative derived functors. This serves as an example to support the assertion of Holm \cite{Holm} that ``every result in classical homological algebra has a counterpart in Gorenstein homological algebra''.

\begin{cor}\label{cor 5}
If $R$ is a left-Gorenstein ring, then $\mathcal{W}$-$\mathrm{g.pd}(R\text{-}\mathrm{Mod})$= $\mathcal{W}$-$\mathrm{g.id}(R\text{-}\mathrm{Mod})$. Moreover, for any ring $R$ we have $$\mathrm{sup}\{\mathrm{G}\text{-}\mathrm{proj.dim}_{R}(M) \mid M\in R\text{-}\mathrm{Mod} \} = \mathrm{sup}\{\mathrm{G}\text{-}\mathrm{inj.dim}_{R}(M) \mid M\in R\text{-}\mathrm{Mod} \}.$$
\end{cor}

When $R$ is left-Gorenstein, we can lift the cotorsion triple $(\mathcal{GP}, \mathcal{W}, \mathcal{GI})$ to $R$-complexes consisting of Gorenstein projective and injective modules. Then we get the following homotopy equivalences.

\begin{cor}\label{cor 6}
Let $R$ be a left-Gorenstein ring. Then there is an equivalence of categories $\mathrm{K}(\mathcal{GP})\simeq \mathrm{K}(\mathcal{GI})$, which restricts to an equivalence $\mathrm{K}(\mathcal{P})\simeq \mathrm{K}(\mathcal{I})$.
\end{cor}

\begin{proof}
The equivalence $\mathrm{K}(\mathcal{GP})\simeq \mathrm{K}(\mathcal{GI})$ follows directly from Theorem \ref{thm 4}. It remains to prove the restricted equivalences.

Specifically, the equivalent functor is $\mathrm{L}:\mathrm{K}(\mathcal{GP})\rightarrow \mathrm{K}(\mathcal{GI})$, which takes fibrant replacement, i.e. special left $dw\widetilde{\mathcal{GI}}$-approximation by using the completeness of $(\widehat{\mathcal{W}}, dw\widetilde{\mathcal{GI}})$. For any complex $M\in dw\widetilde{\mathcal{P}}$, we have an exact sequence $0\rightarrow M\rightarrow G\rightarrow W\rightarrow 0$, where $G\in dw\widetilde{\mathcal{GI}}$ and $W\in\widehat{\mathcal{W}}$. Since $R$ is left-Gorenstein, every projective module has finite injective dimension. Then $M$ and $W$ are complexes with each item having finite injective dimension. Therefore, $G$ is a complex composed by modules with finite injective dimension. Since the injective dimension of any Gorenstein injective module is either zero or infinite \cite[Proposition 10.1.2]{EJ00}, it follows that $G = \mathrm{L}(M)$ is in $dw\widetilde{\mathcal{I}}$. Dually, consider the equivalent functor $\mathrm{R}: \mathrm{K}(\mathcal{GI})\rightarrow \mathrm{K}(\mathcal{GP})$ of taking cofibrant replacement, i.e. special right $dw\widetilde{\mathcal{GP}}$-approximation by using the completeness of $(dw\widetilde{\mathcal{GP}}, \widehat{\mathcal{W}})$. It maps the complex in $dw\widetilde{\mathcal{I}}$ to $dw\widetilde{\mathcal{P}}$. Consequently, we have a restricted equivalence $\mathrm{K}(\mathcal{P})\simeq \mathrm{K}(\mathcal{I})$.
\end{proof}

Denote by $R$-mod the category of finitely presented $R$-modules, and by $\mathrm{D}^{b}(R\text{-}\mathrm{mod})$ its bounded derived category. Let $R^{op}$ be the opposite ring of a ring $R$. We use $\mathrm{K}(\mathcal{P})^{c}$ and  $\mathrm{K}(\mathcal{I})^{c}$ to denote the subcategories consisting of compact objects. Note that there are equivalences $\mathrm{K}(\mathcal{P})^{c}\simeq \mathrm{D}^{b}(R^{op}\text{-}\mathrm{mod})$ by Neeman \cite[Proposition 7.12]{Nee08}; compare to J{\o}rgensen \cite[Theorem 3.2]{Jor05}; and $\mathrm{K}(\mathcal{I})^{c}\simeq \mathrm{D}^{b}(R\text{-}\mathrm{mod})$ by Krause \cite[Proposition 2.3(2)]{Kra05}. Then we get \cite[Theorem C]{Chen10} as follows.

\begin{cor}\label{cor 7}
Let $R$ be a left-Gorenstein ring which is left noetherian and right coherent. Then there is a duality $\mathrm{D}^{b}(R^{op}\text{-}\mathrm{mod})\simeq \mathrm{D}(R\text{-}\mathrm{mod})$ of triangulated categories.
\end{cor}

\begin{rem*}
\indent $(1)$ Consider the cotorsion triple $(\mathcal{P}, R\text{-}\mathrm{Mod}, \mathcal{I})$. By Theorem \ref{thm 4}, if $R$ is a ring of finite global dimension, then we have a trivial equivalence: $\mathrm{K}(\mathcal{P})\simeq \mathrm{D}(R)\simeq \mathrm{K}(\mathcal{I})$.\\
\indent $(2)$ It should be of particular interest to note that over a general ring $R$, a complex $M$ is in $dw\widetilde{\mathcal{GP}}$ if and only if it is a Gorenstein projective graded $R[x]/(x^2)$-modules \cite{Ren18}; if and only if it is a categorically Gorenstein projective complex \cite{Ren18, Y11}, that is, there exists an exact sequence $\mathbb{P}$ of projective complexes, such that $\mathbb{P}$ remains exact after applying $\mathrm{Hom}_{\mathrm{Ch}(R)}(-, P)$ for any projective complex $P$, and $M$ is a cycle of $\mathbb{P}$.
\end{rem*}

\begin{ack*}
The author is grateful to the referee for several suggestions that improved the paper, and he thanks Professor Liu Zhongkui for helpful comments on an early version of the manuscript. This work was supported by National Natural Science Foundation of China (No. 11871125), Natural Science Foundation of Chongqing (No. cstc2018jcyjAX0541) and the Science and Technology Research Program of Chongqing Municipal Education Commission (No. KJQN201800509).\\
\end{ack*}

\end{document}